\documentclass[12pt]{article}
\usepackage{a4wide}

\usepackage{amsfonts}
\usepackage{amscd,amsmath}
\usepackage{euscript}
\newcommand{\doublerightarrow}{\mbox{\Ar\symbol{16}}}

\usepackage[arrow,matrix,curve]{xy}\SilentMatrices
\def\xyma{\xymatrix@M.1em}

\usepackage{epsfig}
\usepackage{amssymb}
\newfont{\Ar}{msam10}
\newcommand{\MySquare}{\hfill{\mbox{\Ar\symbol{3}}}}
\usepackage{psfrag}
\newtheorem{theorem}{Theorem}
\newtheorem{proposition}[theorem]{Proposition}

\newtheorem{definition}[theorem]{Definition}
\newtheorem{conjecture}[theorem]{Conjecture}
\input xy
\xyoption{all}

\begin{document}

\title{Persistent homology of groups}

\author{Graham Ellis \& Simon King\footnote{
Supported by  Marie Curie fellowship MTKD-CT-2006-042685}
\\
{\small School of Mathematics, National University of Ireland, Galway}}
\maketitle

\maketitle

\date{}

\abstract{ 
We introduce and investigate
 notions of persistent homology for $p$-groups and for coclass trees of $p$-groups.     
 Using computer techniques we show that  
 persistent homology provides     
  fairly strong homological   
invariants for   
  $p$-groups of  order $\le 81$.  
The 
 strength of these invariants, and 
 some elementary theoretical properties, 
  suggest 
 that  persistent homology may
 be a useful tool in the study 
of prime-power groups.
  
}

\begin{section}{Introduction}
Persistent homology is a  tool from applied topology
that was introduced for studying
qualitative properties of
 large empirical 
data sets \cite{edelsbrunner}. 
At its simplest, 
the idea 
is 
to impose some  metric  on a data set $X$, choose an
 appropriate  sequence of  
 metric space inclusions
 $X=X_0 \subset X_1 \subset X_2 \subset \cdots \subset X_N$, and 
then study the behaviour of the induced homology maps $H_n(X_0,\mathbb F)
\rightarrow H_n(X_1,\mathbb F) \rightarrow H_n(X_2,\mathbb F) \rightarrow \cdots
\rightarrow H_n(X_N,\mathbb F)$ in a given degree $n$. 
The coefficients $\mathbb F$ are typically 
chosen to be a field.  
 Such a sequence of linear maps is then determined, up to isomorphism, by an upper trianguar matrix 
$P_n=(p_{i,j})$ with $p_{i,j}$ the dimension of the  image of the
map $H_n(X_i,\mathbb F) \rightarrow H_n(X_j,\mathbb F)$. In particular 
$p_{i,i}$ is the dimension of $H_n(X_i,\mathbb F)$. The matrix $P_n$ 
contains information on the extent to which
homology $n$-cycles persist through lengths of the induced
sequence. Cycles that persist for a 
significant length  are deemed to be significant and, for appropriately chosen $X_i$, are 
likely to represent some qualitative feature of the initial
data set $X$. Cycles that persist for only a short length are deemed to be less
 significant.

 Persistent homology analysis
can be applied to any data set $X$ for which we have a suitable topology, 
and for which we have a meaningful sequence of topological
 inclusions. 
It provides a concise set of numerical descriptors for  homological
features  
 of (the sequence of inclusions associated to) $X$.  
 In this paper we investigate the potential of applying the idea to groups and
  particularly to finite $p$-groups.    

\medskip
We can view the elements of a  group $G$  as being 
the vertices of a Cayley graph. 
 Furthermore, we can view the Cayley graph as the $1$-skeleton 
of the universal cover $EG$ of a classifying CW-space $BG$. We set $X=BG$
and construct an inclusion $X\hookrightarrow X_1$
 from any surjective group homomorphism
$\phi\colon G\doublerightarrow Q$ with  
 $X_1=BQ$ a classifying space  obtained by
attaching  cells to $BG$. 
  A sequence of inclusions $X\hookrightarrow X_1\hookrightarrow \ldots \hookrightarrow X_N$  
 corresponds to a sequence  $G\doublerightarrow Q_1 \doublerightarrow Q_2 
\doublerightarrow \ldots
\doublerightarrow Q_N$ of
surjective group homomorphisms, or 
equivalently, to an increasing sequence $N_1\le N_2 \le  \cdots \le N_N$
 of normal subgroups of $G$ (where $N_i$ is the kernel of the composite surjection $G\doublerightarrow Q_i$ ).

We  focus  on finite prime-power groups $G$, and on 
the following five normal series in $G$. 
$$\begin{array}{lllllllll}
L_1(G)&=&G, & L_{i+1}(G)&=&[L_i(G),G] &{\it (lower~ central)}\\
L_1^p(G)&=&G, & L_{i+1}^p(G)&=&[L_i^p(G),G](L_i^p(G))^p &{\it (lower~p-central)} \\
D_1(G)&=&G, & D_{i+1}(G)&=&[D_i(G),D_i(G)] &{\it (derived)}\\
Z_0(G)&=&1, & Z_{i+1}(G)&=&{\rm preimage~of~}Z(G/Z_i(G)){\it ~in~}G &{\it (upper~ central)}\\
Z_0^p(G)&=&1, & Z_{i+1}^p(G)&=&{\rm elements ~of~order}\le p{\it ~in~the}\\
&&&&&{\rm preimage~of~}Z(G/Z_i^p(G)){\rm ~in~}G &{\it (upper~ p-central)}\\
\end{array}$$
These five series can
 be regarded as  functors from the category whose objects are groups and whose arrows are surjections of groups. 
They can be regarded as functors to the category whose objects are 
sequences of group homomorphisms, and whose morphisms are commutative  diagrams of groups. So, for instance, we view $Z$ as a functor
 which sends a surjection $G\rightarrow Q$ to the following
commutative diagram.
$$\begin{array}{cccccccc}
G &\rightarrow &G/Z_1(G) &\rightarrow &G/Z_2(G) & \rightarrow &\cdots\\
\downarrow &&\downarrow&&\downarrow\\
Q &\rightarrow &Q/Z_1(Q) &\rightarrow &Q/Z_2(Q) & \rightarrow &\cdots
\end{array}$$

For $F$ equal to any
 of $L$, $L^p$, $D$, $Z$, $Z^p$ we define the {\it persistence matrix} 
$P_n^F(G) =(p_{i,j})$ to be an upper triangular matrix. For $F$ equal to
 $Z$ or $Z^p$ and $i\ge j$ the   entry
 $p_{i,j}$  is   
the dimension of the  image of the
map $$h_{i,j}\colon H_n(B(G/F_i(G)),\mathbb F_p)  \rightarrow H_n(B(G/F_j(G)),\mathbb F_p).$$
For $F$ equal to  $L$, $L^p$ or $D$ and $i\ge j$ the   entry
 $p_{i,j}$  is
the dimension of the  image of the
map $$h_{i,j}\colon H_n(B(G/F_j(G)),\mathbb F_p)  \rightarrow H_n(B(G/F_i(G)),\mathbb F_p).$$
 The family $H_n^F(G) =\{h_{i,j}\}_{i\ge j}$ is called a  
  {\it persistence module} and is a functorial invariant of the group $G$.
Two persistence modules are isomorphic if and only if the 
 corresponding  persistence matrices are identical. 

  Our aim is to investigate  the extent to which  
 persistence matrices  
 can be used to determine the structure of  finite $p$-groups. 
 For instance, we use computer techniques \cite{kingsage,greenking,hap,ellissmith}
 to establish that the degree seven
 upper central series persistence matrix 
 $P_7^{Z}(G)$ yields 181 distinct matrices when 
$G$ ranges over the 267 groups of order $64$. Furthermore, the groups of order 64 give rise to 187 distinct infinite
 sequences $P_\ast^{Z}(G) = (P_n^Z(G))_{n\ge 1}$.
 We give analogous statistics for each of the five series $F$ and all prime-power groups of order at most $81$.
  We also give some elementary
theoretical results aimed at
  understanding the nature 
 of the group-theoretic information contained in persistence matrices. We believe 
 that the
apparent strength of persistence matrices as group  invariants, and
 their basic theoretical properties,
  suggest
 that  persistent homology may
 be a useful tool in the study
of prime-power groups.

\end{section}

\begin{section}{Examples and properties of persistence}
 Consider the dihedral group $D_{32}$ of order 64. 
Using algorithms recently implemented in the 
group cohomology package \cite{kingsage}  for {\sc sage} (see \cite{greenking} for an overview of  its algorithms)   or the  
{\sc gap} homological algebra package {\sc hap}
\cite{hap}   (see \cite{ellissmith} for an overview
of its algorithms)
one can compute the   lower central series persistence matrix of $D_{32}$ in  degree 2 
 to be 
$$P_2^{L}(D_{32}) = \left( \begin{array}{lllllllll}
 3  &2  &2  &2  &2\\ 
 0  &3  &2  &2  &2\\ 
 0  &0  &3  &2  &2\\ 
 0  &0  &0  &3  &2\\
 0  &0  &0  &0  &3 
\end{array}\right) .$$
The first row of this matrix  implies that $H_2(D_{32},\mathbb F_2)$ has 
dimension $3$, and that precisely two basis elements persist ({\it i.e.} remain non-zero) under the induced maps
 $H_2(D_{32},\mathbb F_2) \rightarrow H_2(D_{32}/L_i(D_{32}),\mathbb F_2)$ ($2\le i\le 5$). 
 The second row implies that
 $H_2(D_{32}/L_5(D_{32}),\mathbb F_2)$ has a basis of three elements, precisely two of which persist
 under the induced maps
 $H_2(D_{32}/L_5(D_{32}),\mathbb F_2) \rightarrow H_2(D_{32}/L_i(D_{32}),\mathbb F_2)$ ($2\le i\le 4$). 
 The matrix is represented by  
 the  {\it persistence  bar code} shown in Figure \ref{PFLD32}(b) and,
in fact, can
 be reconstructed from the information in this bar code. Persistence bar codes
for the matrices   
$P_n^{L}(D_{32})$,  $n=1,3$, 
are given in Figures \ref{PFLD32}(a) and (c). 
\begin{figure}
$${\psfrag{1}{(a)}\psfrag{2}{(b)}\psfrag{3}{(c)}\epsfig{file=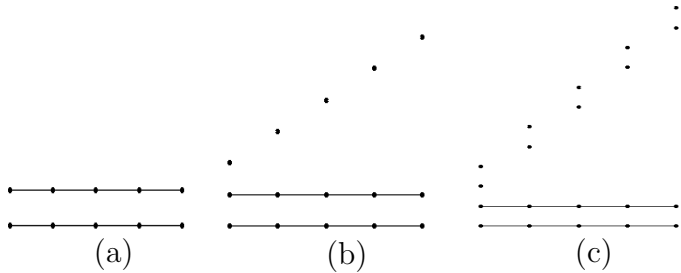,height=3.5cm}}$$
\caption{Degree 1, 2 and 3 lower central bar codes for $D_{32}$ } \label{PFLD32}.
\end{figure}
 Persistence bar codes
for the matrices
$P_n^{L}(Q_{32})$,  $n=1,2,3$, associated to the quaternion group of order 64 
are given in Figure \ref{PFLQ32}.
\begin{figure}
$${\psfrag{1}{(a)}\psfrag{2}{(b)}\psfrag{3}{(c)}\epsfig{file=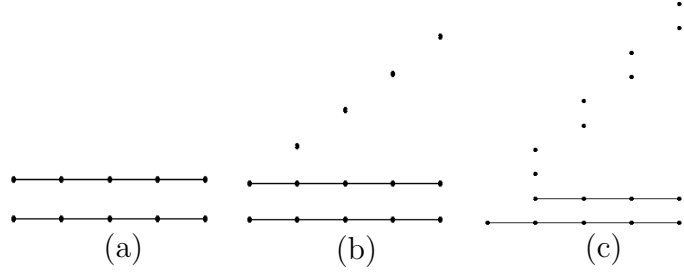,height=3.5cm}}$$
\caption{Degree 1, 2 and 3 lower central bar codes for $Q_{32}$ } \label{PFLQ32}.
\end{figure}
 Persistence bar codes
for the matrices
$P_n^{L}(QD_{32})$,  $n=1,2,3$, associated to the quasi-dihedral
 group of order 64
are given in Figure \ref{PFLQD32}.
 The use of  bar codes for describing persistence matrices
 was  introduced by G. Carlsson {\it et al.} in \cite{carlsson}.
\begin{figure}
$${\psfrag{1}{(a)}\psfrag{2}{(b)}\psfrag{3}{(c)}\epsfig{file=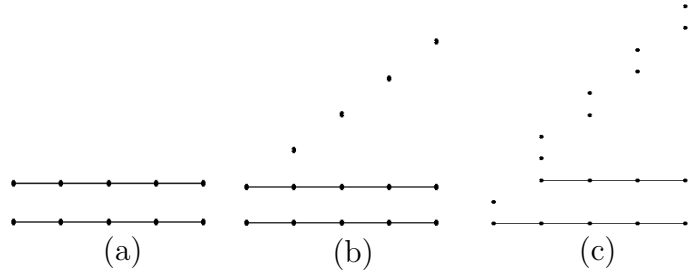,height=3.5cm}}$$
\caption{Degree 1, 2 and 3 lower central bar codes for $QD_{32}$ } \label{PFLQD32}.
\end{figure}

We are interested in the extent to which  group-theoretic  information is
 retained by persistence matrices for the functors
  $L$, $L^p$, $D$, $Z$, $Z^p$. 
  We shall let $F$ denote any
one of the above five functors and let
  $F_i(G)$ denote the $i$th term in the corresponding normal subseries of $G$. 

\begin{proposition}
Let $G$ be  a $p$-group.\newline
(i) For $F$ any of the above five functors the first persistence matrix
$P_1^{F}(G)$ determines 
the minimal number of generators of $G/F_i(G)$ for 
all $i\ge 1$.\newline
(ii) For $F=$ $L$ or $Z$ the first persistence matrix $P_1^{F}(G)$ 
determines the nilpotency class of $G$.\newline
(iii) For $F=$ $L^p$ or $Z^p$ the first two
persistence matrices $P_1^{F}(G)$ and  $P_2^{F}(G)$
determine the order of $G$.
\end{proposition}
{\bf Proof.} Assertion (i) follows from the fact that the 
dimension of the vector space $H_1(G/F_i(G),\mathbb F_p)\cong G/F_i(G)[G,G]G^p$ is equal to the minimal number of generators of $G/F_i(G)$.
 This dimension is the entry $p_{i,i}$ in the first persistence matrix. Assertion (ii) is just the observation that the number of columns in the persistence matrix is by definition equal to  the  length of the upper or lower central series of $G$.
 We prove assertion (iii) just for the functor $F=Z^p$. We use induction 
on the length $k$ of the upper $p$-central series. If  
$k$=1 then $G$ is the trivial group. If $k=2$ then $G$ is an elementary abelian $p$-group of order $p^d$ where $d$
   can be determined by (i).
As an inductive hypothesis suppose that the assertion is true when the upper
 $p$-central series has length $k$. For $G$  a group with upper $p$-central series of length $k+1$  we set $Q=G/Z_1^p(G)$.   The five term natural exact homology sequence
$$H_2(G,\mathbb F_p) \rightarrow H_2(Q,\mathbb F_p) \rightarrow Z_1^p(G) \rightarrow
H_1(G,\mathbb F_p) \rightarrow H_1(Q,\mathbb F_p) \rightarrow 0 \eqno(1)$$ 
allows us to determine the dimension of the vector 
space $Z_1^p(G)$ from the first two upper $p$-central persistence matrices. 
By the inductive hypothesis we can determine the order of $Q$ from 
these two matrices. Then we have the order $|G|=|Q||Z_1^p|$ as required. 
  \MySquare

\begin{proposition}\label{abelian}
The abelian invariants of an abelian $p$-group $G$ are determined by
the first upper $p$-central persistence
 matrix $P_1^{Z^p}(G)$. 
\end{proposition}
{\bf Proof.} We can work by induction on the length $k$
of the upper $p$-central series. If $k=1$ then the persistence matrix
 has just one entry, namely the dimension of the elementary abelian group $G$.
 In general we set $Q=G/Z_1^p(G)$ and, as an inductive hypothesis, assume 
the proposition true for $Q$. Any surjection $G\rightarrow Q$ of abelian groups
induces a surjection in second homology $H_2(G,\mathbb F_p)\rightarrow H_2(Q,\mathbb F_p)$. The exact sequence (1) thus allows
 us to determine the dimension $d$ of the vector space 
$Z_1^p(G)$ from $P_1^{Z^p}(G)$. The abelian invariants of $G$ can be obtained
from those of $Q$  by 
multiplying  precisely $d$ of the highest abelian  invariants of $Q$ by $p$.
 As an example, the bar code for $P_1^{Z^2}(G)$ is given in Figure \ref{abelianbarcode} for the group $G=C_2\times C_4\times C_4\times C_{16}$.
 \MySquare
\begin{figure}
\begin{center}
\epsfig{file=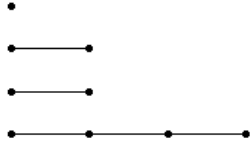,height=2cm}
\caption{First upper 2-central bar code for $G=C_2\times C_4\times C_4\times 
C_{16}$} \label{abelianbarcode}
\end{center}
\end{figure}

\medskip
We regard a bar code as a directed graph whose 
 vertices are arranged in columns, whose   
   edges  connect certain pairs of vertices in neighbouring columns, and where
 all edges are horizontal and directed away from the first column.

\begin{proposition}\label{secondbarcode}
Let $G$ be a finite $p$-group.\newline
(i) The bar code for the first  lower central
persistence matrix $P_1^{L}(G)$ consists of $d={\rm rank}(G/[G,G]G^p)$ horizontal paths, each starting in the first column and ending in the final column.
\newline
(ii) In the bar code for the second  lower central
persistence matrix $P_2^{L}(G)$
every horizontal path  starts in the first column.
\newline
(iii) In the bar code for the second  lower central
persistence matrix $P_2^{L}(G)$
 let $v$ be the number of vertices in the $j$th column ($j\ge 2$) 
 that are not incident with an edge. Set $j'=c+2-j$ where $c$ is the nilpotency class of $G$. 
 Then $v$ is the dimension of the vector space $L_{j'}(G)/L_{j'+1}(G)\otimes \mathbb F_p$.
 (In particular, the number of vertices  in the right-most column 
 not incident with an edge is equal to the rank of $L_2(G)/L_3(G)\otimes \mathbb F_p$.)
\end{proposition}
{\bf Proof.} Assertion (i) follows from the canonical isomorphisms
$H_1(G/L_i(G),\mathbb F_p)  \cong G/[G,G]G^p $ for $i\ge 0$.
Assertion (ii) follows from the  natural exact sequences
$$H_2(G/L_{i+1}(G),\mathbb F_p) \rightarrow H_2(G/L_i(G),\mathbb F_p) \rightarrow L_i(G)/L_{i+1} \otimes \mathbb F_p \rightarrow 0$$
$$H_2(G,\mathbb F_p) \rightarrow H_2(G/L_i(G),\mathbb F_p) \rightarrow L_i(G)/L_{i+1} \otimes \mathbb F_p \rightarrow 0$$
that arise as part of the  five term  exact sequence in mod-$p$ homology.
  These two  sequences imply that  the   two homomorphisms
$H_2(G/L_{i+1}(G),\mathbb F_p) \rightarrow H_2(G/L_i(G),\mathbb F_p)$ and
$H_2(G,\mathbb F_p) \rightarrow H_2(G/L_i(G),\mathbb F_p)$ have identical 
images. This identity implies (ii).
 Asertion (iii) follows from the second of the exact sequences.  \MySquare


\end{section}

\begin{section}{Persistence matrices as  group invariants}

 Proposition \ref{abelian} states that $P_1^{Z^p}$ is a complete
 invariant for abelian $p$-groups. We now investigate the strength of
 persistent homology as a group invariant for finite prime-power
 groups of low order.  The computations were made using the second
 author's group cohomology package~\cite{kingsage} for the computer algebra system Sage~\cite{Sage}.
 Where possible, the computations were corroborated using  the first author's {\sc gap} package \cite{hap}). 
 We begin with the following  summary of the computations for 
$P_\ast^F(G) = (P_n^F(G))_{n\ge 1}$.

\begin{theorem}
For the 366 groups of order at most 81:\newline
(i)   $P_*^{Z}$  
  partitions the  groups  into 277 classes with maximum class size equal to 7.
\newline
(ii) 
   $P_*^{Z^p}$ 
  partitions the  groups  into 262 classes with maximum class size equal to 8.
\newline
(iii)
  $P_*^{L}$
  partitions the  groups  into 227 classes with maximum class size equal to 7.
\newline
(iv)
   $P_*^{L^p}$ 
  partitions the  groups  into 179 classes with maximum class size equal to 15.
\newline
(v)
  $P_*^{D}$ 
  partitions the  groups  into 180 classes with maximum class size equal to 13.
\end{theorem}

A more detailed description of the computations is given in 
Tables \ref{TABLE} and
\ref{TABLE2} 
 which contain, for  prime-powers $k=8,16,27,32,64, 81$, and for each of the five functors $F$:
\begin{enumerate}
\item the number ${\rm Nr}(k)$ of isomorphism classes of groups of order $k$.
\item   the integer pair $(|C|,\max)$ where $|C|$ is
the number of classes of groups of order $k$ distinguished by the invariant $P_\ast^F$,
 and $\max$ is
 the maximum size of a class.
\item the smallest integer $t$ for which the invariant $P_{\ast\le t}^F= 
 (P_n^{F})_{n\le t}$ is as strong as $P_\ast^F$ on the groups of order $k$.
\item an integer triple $(|C|,\max, d)$ where $|C|$ is
the number of classes of groups of order $k$ distinguished by the matrix  
$P_d^F$,
 and $\max$ is
 the maximum size of a class (for some choice of $d$). 

\end{enumerate}

\begin{table}
$$\begin{array}{|l||l|l|l|l|l|l|l|l}
F& &k=8 &k=16 &k=27 &k=32 &k=64
&k=81\\
\hline &&&&&&&\\
&{\rm Nr}(k) &5 &14 &5 &51 &267 &15 \\
&&&&&&&\\
\hline &&&&&&&\\
Z &(|C|,\max) &(5,1) &(13,2) &(5,1)&(44,2) &(187,7) &(14,2)\\
&&&&&&&\\
& t & 3 & 4 & 3 & 5 & 6  & 5\\
&&&&&&&\\
&(|C|,\max,d) &(5,1,3) &(13,2,4) &(5,1,3) &(44,2,7) &(181,7 ,7)&(14,2,5)     \\
&&&&&&&\\
\hline &&&&&&&\\
Z^p &(|C|,\max) &(5,1) &(13,2) &(5,1)&(42,3) &(174,8) &(14,2)\\
&&&&&&&\\
& t & 3 & 4 & 3 & 5 & 6  & 5 \\
&&&&&&&\\
&(|C|,\max,d) &(5,1,3) &(13,2,4) &(5,1,3) &(42,3,5) &(166,8 ,7)&(14,2,5)     \\
&&&&&&&\\
\hline
\end{array}$$
\caption{\label{TABLE}}
\end{table}

\begin{table}
$$\begin{array}{|l||l|l|l|l|l|l|l|l}
F& &k=8 &k=16 &k=27 &k=32 &k=64
&k=81\\
\hline &&&&&&&\\
&{\rm Nr}(k) &5 &14 &5 &51 &267 &15 \\
&&&&&&&\\
\hline &&&&&&&\\
L &(|C|,\max) &(5,1) &(12,2) &(5,1)&(37,3) &(145,7) &(14,2)\\
&&&&&&&\\
& t & 3 & 5 & 3 & 5 & 6  &5 \\
&&&&&&&\\
&(|C|,\max,d) &(5,1,3) &(12,2,4) &(5,1,3) &(37,3,5) &(144,7 ,9)&(14,2,5)     \\
&&&&&&&\\
\hline &&&&&&&\\
L^p &(|C|,\max) &(4,2) &(9,2) &(5,1)&(28,5) &(110,15) &(14,2)\\
&&&&&&&\\
& t & 3 & 4 & 3 & 5 & 6  &5  \\
&&&&&&&\\
&(|C|,\max,d) &(4,2,3) &(9,2,4) &(5,1,3) &(28,5,5) &(109,15 ,9)&(14,2,5)     \\
&&&&&&&\\
\hline
 &&&&&&&\\
D &(|C|,\max) &(5,1) &(10,2) &(5,1)&(29,5) &(108,13) &(14,2)\\
&&&&&&&\\
& t & 3 & 4 & 3 & 5 & 6  & 5\\
&&&&&&&\\
&(|C|,\max,d) &(5,1,3) &(10,2,4) &(5,1,3) &(29,5,7) &(106,13 ,11)&(14,2,5)     \\
&&&&&&&\\
\hline 
\end{array}$$
\caption{\label{TABLE2}}
\end{table}

\begin{table}
\medskip
$$\begin{array}{|l||l|l|l|l|l|l|l|l}
& &k=8 &k=16 &k=27 &k=32 &k=64
&k=81\\
\hline &&&&&&&\\
&{\rm Nr}(k) &5 &14 &5 &51 &267 &15 \\
&&&&&&&\\
\hline &&&&&&&\\
H^{Z^p}_{\ast\le 3} &(|C|,\max) &(5,1) &(13,2) &(5,1)&(46,3) &(188,8) &(14,2)\\
&&&&&&&\\
\hline
\end{array}$$
\caption{\label{table3}}
\end{table}

Persistence bar codes can sometimes distinguish between very similar
groups. Consider for example the two groups $G$, $G'$ of order 64
which are given the identification numbers 158 and 160 in the Small
Groups Library~\cite{BeEiOBr} that is available in {\sc gap}. Their
mod-$2$ cohomology rings $H^\ast(G,\mathbb F_2)$ and
$H^\ast(G',\mathbb F_2)$ have: the same Poincar\'e series; the same
number of generators and relations sorted by degree; the same
depth; the same $a$-invariants. Both groups have nilpotency class $3$ and their upper central
series admit isomorphisms $Z_n(G)\cong Z_n(G')$ ,
$Z_{n+1}(G)/Z_n(G)\cong Z_{n+1}(G')/Z_n(G')$ for $1\le n\le 3$.  Their
$p$-upper central series, lower central series, $p$-lower central
series and derived series admit analogous isomorphisms. However, their
upper central bar codes shown are different in degree $3$ (see Figure
\ref{159160}).
\begin{figure}
$$\epsfig{file=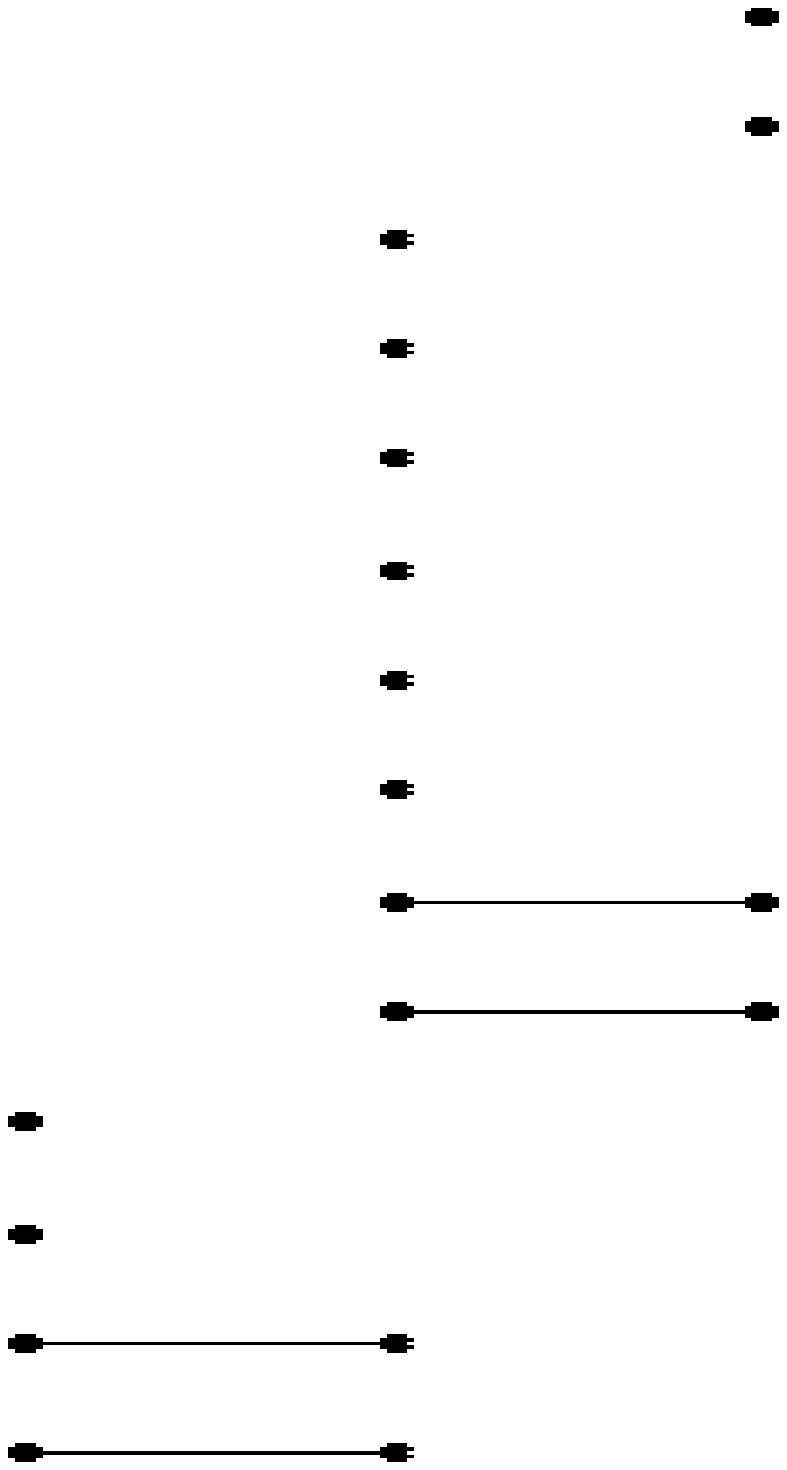,height=3.5cm}$$
\caption{Degree 3 upper central bar codes for groups 158 and 160 of order 64.} 
\label{159160}
\end{figure}

\end{section}

\begin{section}{Integral persistence}

One can use integral homology groups in place of mod $p$ homology  groups when studying persistence. However, the induced homology homomorphisms $H_n(f)\colon
H_n(G,\mathbb Z) \rightarrow H_n(Q,\mathbb Z)$ can not be fully
 described using the notion of dimension. A partial description is
given by the abelian invariants of the source, the target and the cokernel of $H_n(f)$. It is partial due the the extension problem.
 We denote by $IP_n^{F}(G)$ the upper triangular matrix whose entry
in the $i$-th row and $j$-th column ($j\ge i$) is a triple $(A,B,C)$ containing lists of the abelian invariants of the source, target and cokernel of the map
$H_n(G/F_i(G),\mathbb Z) \rightarrow H_n(G/F_j(G),\mathbb Z)$.
  Table \ref{table3} indicates the strength of
$IP^{Z^p}_{\ast\le t}=(IP_n^{Z^p})_{n\le t}$ as an invariant
of the prime-power groups of order at most 81 for $t=3$.

\end{section}
\begin{section}{Persistence in coclass trees}
Recall that a
 $p$-group of order $p^n$ and nilpotency class $c$ is said to have {\it coclass} $r=n-c$. For a fixed $p$ and $r$ one can consider the graph ${\mathbb G}(p,r)$ whose vertices are 
the $p$-groups of coclass $c$. Two groups $G$ and $H$ are connected by an edge in the graph if there exists a normal subgroup $N\le H$ of order $p$
such that $H/N\cong G$. 

The graph  ${\mathbb G}(p,r)$ has infinitely many vertices (since there are infinitely many groups of coclass $r$) and is a forest of trees (since the above $N$ is the smallest non-trivial term of the lower central series of $H$).
 The graph ${\mathbb G}(p,r)$ can be stratified into levels by deeming all groups of order $p^l$ to be at level $l$. 
 The  graph ${\mathbb G}(2,1)$ is shown in Figure
\ref{coclass}. Its three columns contain, respectively, the dihedral groups of order $2^l$, the quaternion groups of order $2^l$ and the semi-dihedral groups of order $2^l$.
\begin{figure}
\begin{center}
\epsfig{file=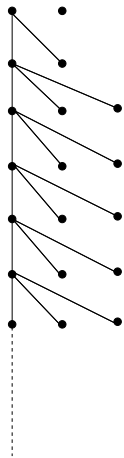,height=6cm}
\caption{Coclass graph $\mathbb G(2,1)$} \label{coclass}
\end{center}
\end{figure}
  Lower central bar codes 
for the three coclass $1$ groups of order 128 
 are shown in Figures \ref{PFLD32}, \ref{PFLQ32}, \ref{PFLQD32}.

Much is known about the graph ${\mathbb G}(p,r)$. A good general reference is the book by Leedham-Green and McKay \cite{leedhamgreenmckay}. It is known
 that the graph is always a forest containing finitely many trees together with
 finitely many sporadic groups. Furthermore, each  tree has
 a unique maximal path of infinite length. 
In the case of $2$-groups it is known that ${\mathbb G}(p,r)$ has 
bounded width ({\it i.e.} there exists some integer $f$ such 
that any vertex is at most $f$ edges away from the infinite 
maximal path). In the case of ${\mathbb G}(2,1)$ the 
infinite path runs through all the dihedral $2$-groups, the sporadic group is the cyclic group of order 4, and the width is $f=2$.

Given a  coclass tree $\mathbb T$ in ${\mathbb G}(p,r)$ we denote by
$S_{\mathbb T}$ the inverse limit of the infinite path. It is
 known that $S_{\mathbb T}$ is a $p$-adic space group. Its
 translation subgroup is an abelian normal subgroup $T \le S_{\mathbb T}$. One
 defines the relative lower central 
series $T_n$ by $T_1=T$ and $T_{n}=[T_{n-1},S_{\mathbb T}]$. 
The quotients $S_{\mathbb T}/T_n$ ($n\ge 1$) are precisely the groups on the infinite path in $\mathbb T$.

 We want to define the persistent homology of a coclass tree $\mathbb T$.     
 Let $G_l$ denote the $p$-group at level $l$ on 
 the  infinite path of a coclass tree $\mathbb T$. 
 Let
${\rm Im}(\nu^{l,k}_n)$  denote the image 
   of the canonical homology homomorphism
$\nu^{l,k}_n\colon H_n(G_{l+k},\mathbb F_p)  \rightarrow H_n(G_l,\mathbb F_p)$. 
\begin{definition}{\rm  
The {\it $l$-persistent homology} of $\mathbb T$ in degree $n$
 is  the subgroup
$$P_lH_n(\mathbb T) = \bigcap_{k=1}^\infty {\rm Im}(\nu^{l,k}_n)$$
of the  degree $n$ homology group $H_n(G_l,\mathbb F_p)$.
}
\end{definition}
Note that 
there is a canonical infinite sequence of surjective homomorphisms
$$\cdots \rightarrow P_{l+2}H_n(\mathbb T) \rightarrow P_{l+1}H_n(\mathbb T)
\rightarrow P_lH_n(\mathbb T)
\ . \eqno(2)  $$

\begin{definition}{\rm 
We define the {\it persistent homology} $PH_n(\mathbb T)$
 of a coclass tree $\mathbb T$ to be the inverse limit
 of the sequence (2) of surjections.
}
\end{definition}

The philosophy is  that $PH_n(\mathbb T)$ should capture some group-theoretic properties that are  common to all groups in the tree.  Easy 
 results
 in this direction are parts (ii) and (iv) of the following proposition. 
Part (i) of the  proposition implies that the surjections in (2)
are isomorphisms for all sufficiently large $l$. 
Hence  $PH_n(\mathbb T)=
{\rm Image}(H_n(G_{l+1},\mathbb F_p) \rightarrow H_n(G_l,\mathbb F_p))$
for all groups $G_l$   on the infinite path 
 in  the tree above a certain level.   
\begin{proposition} (i) The  persistent homology $PH_n(\mathbb T)$ is a finite dimensional vector space for all degrees $n\ge 1$.\newline
(ii) The dimension of $PH_1(\mathbb T)$ equals the minimum number of generators for any group in the tree.
\newline
(iii) 
 $H_2(G,\mathbb F_p) \cong PH_2(\mathbb T) \oplus \mathbb F_p$ for all groups $G$ above a certain level 
  in the tree which are not leaves. For leaves, the dimension of $H_2(G,\mathbb F_p)$ is at least that of $PH_2(\mathbb T)$.
\newline
(iv) 
 Any group $G$ in the tree needs at least $\dim(PH_2(\mathbb T))$ relators to define it. If the group is not a leaf it needs at least $\dim(PH_2(\mathbb T))+1$
relators. 
\end{proposition}
{\bf Proof.}
(i)  The $p$-adic space group associated to 
$\mathbb T$ decomposes into a
short  exact sequence $1\rightarrow T \rightarrow S_{\mathbb T} 
\rightarrow P\rightarrow 1$ where $P$ is a finite (point) group. Each 
group $G_l$
on the infinite path  of $\mathbb T$  thus fits into a short exact sequence
$1\rightarrow T/T_l \rightarrow G_l \rightarrow P\rightarrow 1$. Let $R_\ast^P\rightarrow \mathbb Z$ be any free $\mathbb ZP$-resolution of the integers. Let
 $R_\ast^{(T/T_l)}\rightarrow \mathbb Z$ be the minimal
free $\mathbb Z(T/T_l)$-resolution of the integers constructed as a tensor product of  resolutions of the cyclic summands of $T/T_l$.  Note that the number of
free generators of $R_\ast^{(T/T_l)}$ in a given degree is independent of $l$.
 By a Lemma of C.T.C. Wall the boundary map in the 
tensor product $R_\ast^{(T/T_l)} \otimes R_\ast^P$ can be perturbed to 
produce a free $\mathbb ZG_l$-resolution $R_\ast^{(T/T_l)} \tilde\otimes R_\ast^P$. By construction, the number 
of generators of this latter resolution, in a given degree, 
is independent of $l$. This implies that for a given $n$ the dimension 
of the homology groups $H_n(G_l,\mathbb F_p)$ is bounded by a number depending
only on $P$ and $T$.  This means that the sequence of dimensions  of the 
vectors spaces 
$P_lH_n(\mathbb T), P_{l+1}H_n(\mathbb T), \ldots$ is bounded above. The sequence is also monotonically increasing since the maps in (2) are surjective. Hence the sequence of dimensions converges to the dimension of the inverse limit.

 (ii) This follows directly from the definition of
$PH_1(\mathbb T)$ and the isomorphism $H_1(G,\mathbb F_p) \cong G/[G,G]G^p$. 

(iii)  Let $G_{l+1}\rightarrow G_l$ be a homomorphism in the tree from level $l+1$ to level $l$ with kernel $K$ of order $p$. The five term natural exact homology sequence
$$H_2(G_{l+1},\mathbb F_p) \rightarrow H_2(G_l,\mathbb F_p) \rightarrow K  \rightarrow
H_1(G_{l+1},\mathbb F_p) \stackrel{\cong}{\rightarrow} H_1(G_l,\mathbb F_p) 
\rightarrow 0 \eqno(3)$$
implies 
${\rm Image}(H_n(G_{l+1},\mathbb F_p) \rightarrow H_n(G_l,\mathbb F_p)) \oplus 
K \cong H_2(G_l,\mathbb F_p)$.  If $G_{l+1}$ happens to be on the infinite path in the tree then, by (i),
 the first term in this sum
 stabilizes to $PH_2(\mathbb T)$ for large $l$. If $G_{l+1}$ is not on the infinite path but is not a leaf, then (i) together with Proposition \ref{secondbarcode}(ii) show that
first term in the sum
 stabilizes to $PH_2(\mathbb T)$. If $G$ happens to be a leaf then we can at least conclude from Proposition \ref{secondbarcode} that $H_2(G_{l+1},\mathbb F_p)$
maps onto $PH_2(\mathbb T)$.

(iv) A presentation for $G$  yields the low-dimensional terms in a free $\mathbb ZG$-resolution $R^G_\ast \mathbb Z$ where the $\mathbb ZG$-rank of $R^G_2$ equals the number of relators in the presentation. Clearly this
 rank has to be at least the dimension of $H_2(G,\mathbb F_p)$. So (iii) gives the required result.
 \MySquare

\medskip
The persistent homology could easily be
  computed for some coclass trees. For instance, computer calculations strongly suggest the following result, whose proof should be just a routine
 homological calculation.

\begin{conjecture} For $\mathbb T$ the infinite tree in $\mathbb G(2,1)$ we have 
$$PH_n(\mathbb T) = \mathbb F_2 \oplus \mathbb F_2 \ \ (n\ge 1).$$
\end{conjecture}

\end{section}

\end{document}